\documentclass[12pt]{article}
\usepackage{amssymb,amsthm,amsmath,latexsym}
\usepackage{url}
\usepackage{color}
\usepackage{comment}
\usepackage{lscape}
\newtheorem{thm}{Theorem}
\newtheorem{prop}[thm]{Proposition}
\newtheorem{lem}[thm]{Lemma}

\theoremstyle{remark}
\newtheorem{rem}[thm]{Remark}
\newcommand{\FF}{\mathbb{F}}
\newcommand{\ZZ}{\mathbb{Z}}

\newcommand{\cT}{\mathcal{T}}
\DeclareMathOperator{\wt}{wt}
\DeclareMathOperator{\Aut}{Aut}

\begin{document}
\title{Hadamard matrices of order $36$ formed by codewords
in some ternary self-dual codes}

\author{
Masaaki Harada\thanks{
Research Center for Pure and Applied Mathematics,
Graduate School of Information Sciences,
Tohoku University, Sendai 980--8579, Japan.
email: \texttt{mharada@tohoku.ac.jp}}
and
Keita Ishizuka\thanks{
Research Center for Pure and Applied Mathematics,
Graduate School of Information Sciences,
Tohoku University, Sendai 980--8579, Japan.
email: \texttt{keita.ishizuka.p5@dc.tohoku.ac.jp}}
}

\maketitle

\begin{abstract}
In this note, we study the existence of Hadamard matrices of
order $36$ formed by codewords of weight $36$ in some ternary near-extremal self-dual codes
of length $36$.
\end{abstract}


\section{Introduction}\label{sec:Intro}

Self-dual codes are one of the most interesting classes of codes.
This interest is justified by many combinatorial objects
and algebraic objects related to self-dual codes
(see e.g., \cite{RS-Handbook}).
Among these objects are Hadamard matrices.
In particular, a special class of Hadamard matrices can 
give rise to self-dual codes as their row spaces.
In this note, we  investigate the existence of Hadamard matrices of
order $36$ formed by codewords of weight $36$ in some ternary near-extremal self-dual codes
of length $36$.

A ternary self-dual code $C$ of length $n$
is an $[n,n/2]$ code over the finite field of
order $3$ satisfying
$C=C^\perp$, where $C^\perp$ is the dual code  of $C$.
A ternary self-dual code  of length $n$ exists if and only if $n$ is divisible by four.
It was shown in~\cite{MS-bound} that
the minimum weight $d$ of a ternary self-dual code of length $n$
is bounded by $d\leq 3 \lfloor n/12 \rfloor+3$.
If $d=3\lfloor n/12 \rfloor+3$ (resp.\ $d=3\lfloor n/12 \rfloor$), 
then the code is called {extremal} (resp.\ near-extremal).
For length $36$, the Pless symmetry code is a currently known extremal self-dual code.

Recently, Tonchev~\cite{T} studied Hadamard matrices of order $n$
formed by codewords of weight $n$ in ternary extremal self-dual codes of length
$n$. 
The Pless symmetry code contains a Hadamard matrix~\cite[Theorem~4.2]{Pless}.
This matrix is often called the Paley-Hadamard matrix of type~II (see~\cite{T}). 
Tonchev~\cite{T}  showed that the Pless symmetry
code of length $36$ contains exactly two inequivalent Hadamard matrices of
order $36$, namely,
the Paley-Hadamard matrix of type~II  
and a regular Hadamard matrix.
Inspired by his result, 
we study the existence of Hadamard matrices of
order $36$ formed by codewords of weight $36$ in some ternary near-extremal self-dual codes
of length $36$.

This note is organized as follows.
In Section~\ref{Sec:2},
we give definitions and some known results
of ternary self-dual codes and Hadamard matrices. 
In Section~\ref{Sec:36}, we investigate 
ternary self-dual codes generated by the rows of Hadamard matrices
of order $36$ with automorphism of order $17$ classified in~\cite{T86}.
We also give a classification of ternary near-extremal bordered double circulant self-dual codes,
four-negacirculant self-dual codes and quasi-twisted self-dual codes of length $36$.
In Section~\ref{Sec:Results},
we study the existence of Hadamard matrices
formed by codewords of weight $36$ in 
the ternary  near-extremal self-dual codes of length $36$ classified in Section~\ref{Sec:36}.
Consequently,
$89$ inequivalent Hadamard matrices of order $36$ are constructed from
codewords of weight $36$ in the ternary near-extremal self-dual codes.
We also investigate ternary codes  generated by the rows of the $89$ Hadamard matrices
and their transposes.
Finally, we give a relationship between the Pless symmetry code
and a certain near-extremal four-negacirculant self-dual code
via the regular Hadamard matrix given in~\cite[Theorem~2.4 (i)]{T}.

\section{Preliminaries}\label{Sec:2}

In this section, we give definitions and some known results
of ternary self-dual codes and Hadamard matrices used in this note.

\subsection{Ternary self-dual codes}
Let $\FF_3=\{0,1,2\}$ denote the finite field of order $3$.
A \emph{ternary} $[n,k]$ \emph{code} $C$ is a $k$-dimensional vector subspace
of $\FF_3^n$.
All codes in this note are ternary, and the expression ``codes'' means ``ternary codes''.
The parameters $n$ and $k$ are called the \emph{length} and \emph{dimension} of $C$,
respectively.
A \emph{generator matrix} of $C$
is a $k \times n$ matrix whose rows are a basis of $C$.
The \emph{weight} $\wt(x)$ of a vector $x$ of $\FF_3^n$ is 
the number of non-zero components of $x$.
A vector of $C$ is called a \emph{codeword}.
The minimum non-zero weight of all codewords in $C$ is called
the \emph{minimum weight} of $C$.
The polynomial $\sum_{c \in C} y^{\wt(c)}$ is called
the \emph{weight enumerator} of $C$.

The {\em dual} code $C^{\perp}$ of a code
$C$ of length $n$ is defined as
$
C^{\perp}=
\{x \in \FF_3^n \mid x \cdot y = 0 \text{ for all } y \in C\},
$
where $x \cdot y$ is the standard inner product.
A code $C$ is \emph{self-dual} if $C=C^\perp$.
A self-dual code  of length $n$ exists if and only if $n$ is divisible by four.
Two codes $C$ and $C'$ are \emph{equivalent} if there is a
monomial matrix $P$ over $\mathbb{F}_3$ with $C' = C \cdot P$,
where $C \cdot P = \{ x P\mid  x \in C\}$.
We denote two equivalent codes $C$ and $D$ by $C \cong D$.
All self-dual codes were classified in~\cite{CPS}, \cite{HM}, \cite{MPS} 
and~\cite{PSW} for lengths up to $24$.
 
By the Gleason theorem (see~\cite{MS-bound}),
the weight enumerator of an extremal self-dual code is uniquely determined
for each length.
Moreover, it was shown in~\cite{MS-bound} that
the minimum weight $d$ of a self-dual code of length $n$
is bounded by $d\leq 3 \lfloor n/12 \rfloor+3$.
If $d=3\lfloor n/12 \rfloor+3$, then the code is called \emph{extremal}.
For $n \in \{4,8,12,\ldots,64\}$,
it is known that there is an extremal self-dual code of length $n$ (see~\cite[Table~6]{Huffman}).
For length $36$,
the Pless symmetry code $P_{36}$ (see Section~\ref{Sec:P36} for the definition)
is a currently known extremal self-dual code.
If $d=3\lfloor n/12 \rfloor$, then the code is called \emph{near-extremal}.
Recently, some restrictions on the weight enumerators of 
near-extremal self-dual codes of length divisible by $12$ are given in~\cite{AH}.
For example, the possible weight enumerators of 
near-extremal self-dual codes of length $36$ are as follows:
\begin{equation}\label{eq:W}
\begin{split}
&
1
+ \alpha y^9
+( 42840 - 9 \alpha)y^{12}
+( 1400256 + 36 \alpha)y^{15}
\\&
+( 18452280- 84 \alpha)y^{18}
+( 90370368 + 126 \alpha)y^{21}
\\&
+( 162663480 - 126 \alpha)y^{24}
+( 97808480 + 84 \alpha)y^{27}
\\&
+( 16210656 - 36 \alpha)y^{30}
+( 471240 + 9 \alpha)y^{33}
+( 888 - \alpha)y^{36},
\end{split}
\end{equation}
where $\alpha =8\beta$ and
$\beta \in \{1,2,\ldots,111\}$~\cite{AH}.
Note that the number of codewords of weight $9$ 
in a near-extremal self-dual code $C$ of length $36$
determines the weight enumerator of $C$.

\subsection{Hadamard matrices}
A \emph{Hadamard} matrix $H$ of order $n$ is an $n \times n$ matrix
whose entries are from $\{ 1,-1 \}$ such that $H H^T = nI_n$,
where $H^T$ is the transpose of $H$ and $I_n$ is the identity matrix of order $n$.
Hadamard matrices are widely investigated combinatorial matrices (see e.g., \cite{SeberryYamada}). 
Hadamard matrices are also one of  combinatorial matrices  related to self-dual codes.

It is known that
the order $n$ is necessarily $1,2$, or a multiple of $4$.
Two Hadamard matrices $H$ and $K$ are said to be \emph{equivalent}
if there are $(1,-1,0)$-monomial matrices $P$ and $Q$ with $K = PHQ$.
We denote two equivalent Hadamard matrices $H$ and $K$ by $H \cong K$.
An \emph{automorphism} of a Hadamard matrix $H$ is
an equivalence of $H$ to itself.
The set of all automorphisms of $H$ forms a group under
composition called the \emph{automorphism group} $\Aut(H)$ of $H$.
For orders up to $32$, all inequivalent Hadamard matrices are known
(see~\cite{Kimura} and~\cite{KT32} for orders $28$ and $32$, respectively).
Many Hadamard matrices of order $36$ are known (see e.g., \cite{SeberryHttp}).
Tonchev~\cite{T86} showed that there are $11$ inequivalent Hadamard matrices
of order $36$ with automorphism of order $17$.
Note that one of them is the well-known Paley-Hadamard matrix of  type~II 
(see~\cite{T} for the definition).

\subsection{Pless symmetry code $P_{36}$ of length $36$}\label{Sec:P36}

An $n \times n$  matrix of the following form
\[
\left( \begin{array}{cccccc}
r_0&r_1&r_2& \cdots &r_{n-2} &r_{n-1}\\
\mu r_{n-1}&r_0&r_1& \cdots &r_{n-3}&r_{n-2} \\
\mu r_{n-2}&\mu r_{n-1}&r_0& \cdots &r_{n-4}&r_{n-3} \\
\vdots &\vdots & \vdots &&\vdots& \vdots\\
\vdots &\vdots & \vdots &&\vdots& \vdots\\
\mu r_1&\mu r_2&\mu r_3& \cdots&\mu r_{n-1}&r_0
\end{array}
\right)
\]
is called \emph{circulant} and \emph{negacirculant}
if $\mu =1$ and $-1$, respectively.
If $R$ is an $(n-1) \times (n-1)$ circulant matrix, then
a  $[2n,n]$ code with generator matrix of the following form
\begin{equation}\label{eq:bDCC}
\left(\begin{array}{ccccccccc}
{} & {} & {}      & {} & {} &0& 1 & \cdots &1 \\
{} & {} & {}      & {} & {} &1& {} & {}     &{} \\
{} & {} & I_{n} & {} & {} &\vdots& {} & R     &{} \\
{} & {} & {}      & {} & {} &1& {} &{}      &{} \\
\end{array}\right)
\end{equation}
is called a \emph{bordered double circulant} code.
The Pless symmetry code $P_{36}$ is defined as the bordered double circulant code
with generator matrix of the form~\eqref{eq:bDCC}, where
\[
(0,1,1,2,1,2,2,2,1,1,2,2,2,1,2,1,1)
\]
is the first row of $R$~\cite{Pless}.

Now we give other construction methods of the Pless symmetry code $P_{36}$.
If $N$ is an $n \times n$ negacirculant matrix, then
a $[2n,n]$ code with generator matrix of the following form
\begin{equation}\label{eq:QT}
\left(\begin{array}{cc}
I_n & N
\end{array}\right)
\end{equation}
is called a \emph{quasi-twisted} code.
There is a quasi-twisted self-dual code $C$ such that $C \cong P_{36}$~\cite{HHKK}.
If  $A$ and $B$ are $n \times n$ negacirculant matrices, then
a $[4n,2n]$ code with generator matrix of the following form
\begin{equation}\label{eq:4nega}
\left(
\begin{array}{ccc@{}c}
\quad & {I_{2n}} & \quad &
\begin{array}{cc}
A & B \\
2B^T & A^T
\end{array}
\end{array}
\right)
\end{equation}
is called a \emph{four-negacirculant} code.
There is a four-negacirculant self-dual code $C$ such that $C \cong P_{36}$~\cite{HHKK}.

Let $H$ be a Hadamard matrix of order $n$.
Throughout this note,
we denote by $C(H)$ the code generated by the rows of $H$,
where entries of $H$ are regarded as elements of $\FF_3$.
As described above,
Tonchev~\cite{T86} showed that there are $11$ inequivalent Hadamard matrices
of order $36$ with automorphism of order $17$.
We denote the Hadamard matrix corresponding to the $3$-$(36,18, 8)$ design
$E_i$ in~\cite[Table~4]{T86} by $T_{36,i}$ $(i=1,2,\ldots,11)$.
Note that $T_{36,11}$ is equivalent to the Paley-Hadamard matrix of  type~II.
This implies that $C(T_{36,11}) \cong P_{36}$.

Here we regard the Pless symmetry code $P_{36}$ as a member of the following four
families of codes:
\begin{itemize}
\item[(F1)]
codes $C(T_{36,i})$ from the Hadamard matrices  $T_{36,i}$ $(i=1,2,\ldots,11)$,
\item[(F2)]
bordered double circulant self-dual codes,
\item[(F3)]
four-negacirculant self-dual codes,
\item[(F4)]
quasi-twisted  self-dual codes.
\end{itemize}

\section{Near-extremal self-dual codes of length $36$}
\label{Sec:36}

For the above four families (F1)--(F4), we consider near-extremal self-dual codes of length $36$.
All computer calculations 
in the rest of  this note
were done with the help of \textsc{Magma}~\cite{Magma}.

\subsection{Self-dual codes $C(T_{36,i})$ from $T_{36,i}$}

We verified that $C(T_{36,i})$ are self-dual codes of minimum weight $9$
$(i=1,2,\ldots,10)$.
In addition, we verified that  $C(T_{36,2}) \cong C(T_{36,9})$ and there is no other pair
of equivalent codes among $C(T_{36,i})$.
The numbers $A_9$ of codewords of weight $9$ in $C(T_{36,i})$ are listed
in Table~\ref{Tab:T36wt9}.
These numbers determine the weight enumerators of  $C(T_{36,i})$
(see~\eqref{eq:W}).

\begin{table}[thb]
\caption{Numbers $A_9$ of codewords of weight $9$ in $C(T_{36,i})$}
\label{Tab:T36wt9}
\centering
\medskip
{\small
\begin{tabular}{c|l}
\noalign{\hrule height1pt}
$A_9$ & \multicolumn{1}{c}{$i$} \\
\hline
136 &$1,3,5,6,7,8,10$ \\
544 & $2,4,9$ \\
\noalign{\hrule height1pt}
\end{tabular}
}
\end{table}

\subsection{Bordered double circulant self-dual  codes}

By exhaustive search, we found all distinct near-extremal bordered double circulant 
self-dual codes of length $36$, which must be checked further for equivalences.
By equivalence testing, 
we verified that
there are exactly $12$ inequivalent  near-extremal bordered double circulant 
self-dual codes. 
We denote these codes by $D_{36,i}$ ($i=1,2,\ldots,12$).
For the codes $D_{36,i}$,
the first rows $r$ of the circulant matrices $R$ in their generator matrices~\eqref{eq:bDCC}
are listed in Table~\ref{Tab:bDCC}.
The numbers $A_9$ of codewords of weight $9$ in $D_{36,i}$ are listed
in Table~\ref{Tab:bDCCwt9}.
These numbers determine the weight enumerators of  $D_{36,i}$
(see~\eqref{eq:W}).

\begin{table}[thb]
\caption{Bordered double circulant codes $D_{36,i}$}
\label{Tab:bDCC}
\centering
\medskip
{\small
\begin{tabular}{c|c||c|c||c|c}
\noalign{\hrule height1pt}
$i$ &  $r$  & $i$ &  $r$  & $i$ &  $r$  \\
\hline
 1 & 21022212010000011 & 5 & 10210211100100011 & 9 & 01111221220100112 \\
 2 & 11122112102100021 & 6 & 20211112121100012 &10 & 20111111001000012 \\
 3 & 12121111001000010 & 7 & 12121120112100012 &11 & 10011210000000021 \\
 4 & 02212101201000012 & 8 & 12222022210100121 &12 & 01212112121100012 \\
\noalign{\hrule height1pt}
\end{tabular}
}
\end{table}

\begin{table}[thb]
\caption{Numbers $A_9$ of codewords of weight $9$ in $D_{36,i}$}
\label{Tab:bDCCwt9}
\centering
\medskip
{\small
\begin{tabular}{c|l}
\noalign{\hrule height1pt}
$A_9$ & \multicolumn{1}{c}{$i$} \\
\hline
136 &$1,2,\ldots,7$ \\
408 & $8, 9, 10$ \\
544 & $11,12$ \\
\noalign{\hrule height1pt}
\end{tabular}
}
\end{table}

\subsection{Four-negacirculant self-dual codes}

By exhaustive search, we found all distinct near-extremal four-negacirculant self-dual
codes of length $36$, which must be checked further for equivalences.
By equivalence testing, 
we verified that
there are exactly $260$ inequivalent  near-extremal  four-negacirculant  self-dual 
codes.
We denote these codes by $F_{36,i}$ ($i=1,2,\ldots,260$).
For the codes $F_{36,i}$,
the first rows $r_A$ (resp.\ $r_B$) of the negacirculant matrices $A$ (resp.\ $B$) in their generator matrices~\eqref{eq:4nega}
are listed in Tables~\ref{Tab:4nega1} and~\ref{Tab:4nega2}.
The numbers $A_9$ of codewords of weight $9$ in $F_{36,i}$ are listed
in Table~\ref{Tab:4negawt9}.
These numbers determine the weight enumerators of  $F_{36,i}$
(see~\eqref{eq:W}).

\begin{table}[thb]
\caption{Numbers $A_9$ of codewords of weight $9$ in $F_{36,i}$}
\label{Tab:4negawt9}
\centering
\medskip
{\small
\begin{tabular}{c|l||c|l||c|l}
\noalign{\hrule height1pt}
$A_9$ & \multicolumn{1}{c||}{$i$} &
$A_9$ & \multicolumn{1}{c||}{$i$} &
$A_9$ & \multicolumn{1}{c}{$i$}  \\
\hline
 72 &$ 1, 2,\ldots, 6 $  &312 &$ 105, 106,\ldots, 121 $ &576 &$ 246, 247,248, 249 $ \\
 96 &$ 7, 8 $          &360 &$ 122, 123,\ldots, 154 $ &600 &$ 250, 251 $ \\
144 &$ 9, 10,\ldots, 22 $ &384 &$ 155, 156,\ldots, 171 $ &648 &$ 252, 253,254, 255 $ \\
168 &$ 23, 24, 25, 26 $&432 &$ 172, 173,\ldots, 196 $ &720 &$ 256, 257, 258 $ \\
216 &$ 27, 28,\ldots, 45 $&456 &$ 197, 198,\ldots, 213 $ &744 &$ 259, 260 $ \\
240 &$ 46, 47,\ldots, 60 $&504 &$ 214, 215,\ldots, 237 $ & &\\
288 &$ 61, 62,\ldots,104 $&528 &$ 238, 239,\ldots, 245 $ &&\\
\noalign{\hrule height1pt}
\end{tabular}
}
\end{table}

\subsection{Quasi-twisted self-dual  codes}

By exhaustive search, we found all distinct near-extremal quasi-twisted self-dual
codes of length $36$, which must be checked further for equivalences.
By equivalence testing, 
we verified that
there are exactly $55$ inequivalent  near-extremal quasi-twisted self-dual 
codes.
We denote these codes by $Q_{36,i}$ ($i=1,2,\ldots,55$).

\subsection{Equivalent codes}
We investigate equivalences among the near-extremal self-dual
codes $C(T_{36,i})$, $D_{ 36,i}$, $F_{36,i}$ and $Q_{36,i}$.
We verified that
\[
\begin{array}{lll}
D_{36,  1} \cong  C(T_{36, 3}),&
D_{36,  2} \cong  C(T_{36,  8}),&
D_{36,  3} \cong  C(T_{36,  5}),\\
D_{36,  4} \cong  C(T_{36,  6}),&
D_{36,  5} \cong  C(T_{36,  7}),&
D_{36,  6} \cong  C(T_{36,  1}),\\
D_{36,  7} \cong  C(T_{36, 10}),&
D_{36, 11} \cong  C(T_{36,  4}),&
D_{36, 12} \cong  C(T_{36,  2}) \cong C(T_{36, 9}).
\end{array}
\]
For each code $Q_{36,i}$ ($i=1,2,\ldots,55$),
we verified that there is a code $F_{36,j}$ such that $ Q_{36,i} \cong F_{36,j}$.
In addition, we verified that there is no code $F_{36,i}$ such that $F_{36,i} \cong C$
for each code  $C=C(T_{36, j})$ $(j=1,2,\ldots,8,10)$ and $C=D_{36,j}$ $(j=8,9,10)$.
Hence, we consider the following $272$ inequivalent near-extremal self-dual codes:
\[
\begin{array}{ll}
C(T_{ 36, i}) & (i=1,2,\ldots,8,10), \\
D_{ 36,i} & (i=8,9,10), \\
F_{36,i} & (i=1,2,\ldots,260),
\end{array}
\]
in the next section.

\section{Hadamard matrices in some near-extremal self-dual codes}
\label{Sec:Results}


In this section, 
we  study the existence of Hadamard matrices
formed by codewords of weight $36$ in $C \in {\boldsymbol{C}_{36}}$,
where
\begin{equation*}\label{eq:C}
{{\boldsymbol{C}_{36}}}=
\left\{
\begin{array}{l}
C(T_{ 36, 1}),C(T_{ 36, 2}),\ldots,C(T_{ 36, 8}),C(T_{ 36, 10}),\\
\qquad D_{ 36,8},D_{ 36,9},D_{ 36,10},F_{36,1},F_{36,2},\ldots,F_{36,260}
\end{array}
\right\}.
\end{equation*}

\subsection{Method}\label{sec:Method}

The following lemma is trivial.

\begin{lem}\label{lem:orth}
Let $r_1$ and $r_2$ be vectors of $\{1,-1\}^n \subset \ZZ^n$.
Let $n_1(r_i)$ denote the number of $1$'s in $r_i$ $(i=1,2)$.
Suppose that $n \equiv 0 \pmod 4$ and  $r_1 \ne \pm r_2$.
If $r_1$ and $r_2$ are orthogonal, then $n_1(r_1) \equiv n_1(r_2) \pmod 2$.
\end{lem}

We describe how to find all Hadamard matrices 
formed by codewords of weight $36$ in a given
near-extremal self-dual code $C$ of length $36$.
Suppose that $C$ has at least $72$ codewords of weight $36$.
For the set $W(C)$ of codewords of weight $36$ in $C$,
we define the following sets:
\begin{align*}
W_1(C)&=\{x=(1,x_2,x_3,\ldots,x_{36}) \mid x \in W(C)\}, \\
W_{1,i}(C)&=\{ x \in W_1(C) \mid n_1(x) \equiv i \!\!\!\pmod 2\} \ (i=0,1).
\end{align*}
%
For a vector $x=(x_1,x_2,\ldots,x_{36})$ of $\FF_3^{36}$, we define the following:
\begin{equation}\label{eq:rho}
\overline{x}=(\overline{x_1},\overline{x_2},\ldots,\overline{x_{36}}) \in \ZZ^{36},
\end{equation}
where 
$\overline{0}$,
$\overline{1}$ and 
$\overline{2}$ are $0, 1$  and  $-1 \in \ZZ$, respectively.
%
%
%
For a given set $W_{1,i}(C)$ $(i=0,1)$, a simple undirected graph
$\Gamma_i(C)$ is constructed as follows:
\begin{enumerate}
\item the vertex set is 
$\{\overline{x}  \mid x \in W_{1,i}(C)\}  \subset \ZZ^{36}$,
\item two vertices $x$ and $y$ are adjacent
if $x$ and $y$ are orthogonal.
\end{enumerate}
%
%
%
Clearly, a $36$-clique in $\Gamma_i(C)$ gives a Hadamard matrix for $i=0$ and $1$.
By Lemma~\ref{lem:orth}, all Hadamard matrices 
formed by codewords of $W_{1}(C)$ are obtained
by finding all $36$-cliques in $\Gamma_i(C)$ for $i=0$ and $1$.
For any Hadamard matrix $H$ formed by codewords of weight $36$ in $C$,
it follows that there is a Hadamard matrix $H'$ constructed by the above method
satisfying that $H \cong H'$.


\subsection{Results}

By the above method, 
for each code $C \in {\boldsymbol{C}_{36}}$,
we found all distinct Hadamard matrices formed by codewords of weight $36$ in $C$.
The results are given by using the following notations:
\begin{align*}
N_i(C) &= \text{ the number  of $36$-cliques in $\Gamma_i(C)$ $(i=0,1)$}, \\
N(C)   &= \text{ the number of inequivalent  Hadamard matrices in $C$}, \\
\cT(C)  &=(|W_{1,0}(C)|, |W_{1,1}(C)|, N_0(C), N_1(C),N(C)).
\end{align*}


Firstly, for each code $C \in {\boldsymbol{C}_{36,1}}$,
the tuples  $\cT(C)$ are listed in Table~\ref{Tab:C:H},
where
\[
{\boldsymbol{C}_{36,1}}=
\{C(T_{ 36, 1}),C(T_{ 36, 2}),\ldots,C(T_{ 36, 8}),C(T_{ 36, 10}),\\
D_{ 36,8},D_{ 36,9},D_{ 36,10}\}.
\]
Since $N(C(T_{36, i}))=1$ $(i=1,3,4,\ldots,8,10)$,
we denote by $H_{T,i}$ the Hadamard matrix 
formed by codewords of weight $36$ in $C(T_{36, i})$.
It is easy to see that $T_{36, i} = H_{T,i}$ for $i=1,3,4,\ldots,8,10$.
Since $N(C(T_{36, 2}))=2$ and $C(T_{36, 2}) \cong C(T_{36, 9})$,
$T_{36, 2}$ is one of the two inequivalent Hadamard matrices
formed by codewords of weight $36$ in $C(T_{36, 2})$, and
the other Hadamard matrix $H$ satisfies $H \cong T_{36, 9}$.
Therefore, we have the following:
\begin{prop}\label{prop:1}
\begin{enumerate}
\item
There are ten inequivalent Hadamard matrices 
formed by codewords of weight $36$ in the nine codes $C(T_{36, i})$ $(i=1,2,\ldots,8,10)$.
\item
There is no Hadamard matrix
formed by codewords of weight $36$ in the three codes $D_{ 36,i}$ $(i=8,9,10)$.
\end{enumerate}
\end{prop}

\begin{table}[thb]
\caption{Tuples $\cT(C)$ for $C \in {\boldsymbol{C}_{36,1}}$}
\label{Tab:C:H}
\centering
\medskip
{\small
\begin{tabular}{c|c}
\noalign{\hrule height1pt}
$C$ &  $\cT(C)$\\
\hline
$C(T_{36, i})$ $(i=1,3,10)$ & $(308,  68, 1, 0, 1)$\\
$C(T_{36,  2})$                   & $(104,  68, 2, 0, 2)$\\
$C(T_{36,  4})$                   & $(104,  68, 1, 0, 1)$\\
$C(T_{36,  i})$ $(i=5,6,7,8)$&$( 104, 272, 1, 0, 1)$\\
$D_{36, i}$ $(i=8,9,10)$      &$( 204, 36, 0, 0, 0)$\\
\noalign{\hrule height1pt}
\end{tabular}
}
\end{table}

Secondly,
for each code 
$C \in 
{\boldsymbol{C}_{36}} \setminus{\boldsymbol{C}_{36,1}}
$, 
we found all distinct Hadamard matrices 
formed by codewords of weight $36$ in $C$ by the above method.
It holds that
\[
|\{\cT(F_{36,1}),\cT(F_{36,2}),\ldots,\cT(F_{36,260})\}|=153,
\]
where the $153$ distinct tuples $\cT(F_{36,i})$ are listed in  Tables~\ref{Tab:tuple}
and~\ref{Tab:tuple2}.
For each tuple $\cT$ in the tables, we list all values $i$ such that $\cT=\cT(F_{36,i})$.
For each $i=1,2,\ldots,260$,
we denote by $H_{F,i,i_j}$ the Hadamard matrices 
formed by codewords of weight $36$ in the $260$ codes 
$F_{36, i}$, where $i_j=1,2,\ldots,N(F_{36, i})$.
We verified that
\[
\begin{array}{ll}
H_{F,11, 1} \cong H_{F, 67, 1} \cong H_{F,139, 1},&
H_{F,27, 1} \cong H_{F, 83, 1} \cong H_{F,237, 1},\\
H_{F,38, 2} \cong H_{F, 79, 1} \cong H_{F,177, 1},&
H_{F,49, 2} \cong H_{F,168, 1} \cong H_{F,169, 1},\\
H_{F,56, 2} \cong H_{F,115, 1} \cong H_{F,213, 1},&
H_{F,64, 1} \cong H_{F,172, 1} \cong H_{F,174, 1},
\end{array}
\]
and there is no other pair of equivalent Hadamard matrices among the Hadamard
matrices $H_{F,i,i_j}$.
Therefore, we have the following:

\begin{prop}\label{prop:2}
There are $89$ inequivalent Hadamard matrices 
formed by codewords of weight $36$ in the $260$ codes 
$F_{36,i}$  $(i=1,2,\ldots,260)$.
\end{prop}

\begin{rem}
The above $89$ inequivalent Hadamard matrices can be obtained electronically from
\url{https://www.math.is.tohoku.ac.jp/~mharada/F3-36/Hmat-F.txt}.
\end{rem}

We consider equivalences among the Hadamard matrices
formed by codewords of weight $36$ in the $272$ codes of ${\boldsymbol{C}_{36}}$.
The orders of the automorphism groups $\Aut(H_{F,i,i_j})$ of $H_{F,i,i_j}$ are listed in
Table~\ref{Tab:Aut}.
From the table,
there is no pair of inequivalent Hadamard matrices among 
the $99$ inequivalent Hadamard matrices given in 
Propositions~\ref{prop:1} and~\ref{prop:2}.

\begin{table}[thb]
\caption{Orders of the automorphism groups $\Aut(H_{F,i,i_j})$}
\label{Tab:Aut}
\centering
\medskip
{\small
\begin{tabular}{c|l}
\noalign{\hrule height1pt}
$|\Aut(H_{F,i,i_j})|$ &  \multicolumn{1}{c}{$(i,i_j)$}\\
\hline
8 &$(8,1),(18,2),(23,1),(23,2),(168,2),(168,3)$\\
12&$(7,1),(39,1)$\\
16&$(37,1),(37,3)$\\
24&$(8,2),(18,3)$\\
36&$(4,1),(4,2),(9,1),(11,1),(12,1),(12,2),(13,1),(18,1),(26,1),$\\&
$(26,2),(27,1),(28,1),(30,1),(35,1),(38,1),(38,2),(42,1),$\\&
$(53,1),(56,1),(56,2),(60,1),(61,1),(66,1),(73,1),(96,1),$\\&
$(103,1),(104,1),(110,1),(111,1),(112,1),(117,1),(120,1),$\\&
$(130,1),(133,1),(148,1),(149,1),(166,1),(173,1),(176,1),$\\&
$(176,2),(178,1),(193,1),(197,1),(201,1),(208,1),(210,1),$\\&
$(213,2),(226,1),(227,1),(230,1),(234,1),(235,1),(245,1),$\\&
$(246,1),(247,1),(254,1)$\\
72&$(23,3),(49,1),(49,2),(54,1),(54,2),(64,1),(68,1),(90,1),$\\&
$(108,1),(108,2),(143,1),(143,2),(164,1),(168,4),(184,1),$\\&
$(188,1),(215,1),(224,1),(260,1)$\\
144&$(37,2),(194,1)$\\
\noalign{\hrule height1pt}
\end{tabular}
}
\end{table}

Now we consider equivalences among the $89$ inequivalent Hadamard matrices
$H_{F,i,j}$ and their transposes.
We verified that $H_{F,s,t}^T$ is equivalent to none of 
the $89$ inequivalent Hadamard matrices
$H_{F,i,j}$ for each $(s,t) \in \mathcal{S}_{s,t}$
and there is a Hadamard matrix $H_{F,i,j}$ such that $H_{F,i,j} \cong H_{F,s,t}^T$
for each $(s,t) \not\in \mathcal{S}_{s,t}$, where
\begin{equation}\label{eq:T}
\mathcal{S}_{s,t}=
\left\{
\begin{array}{l}
(7,1),(8,1),(8,2),(18,2),(18,3),(23,1),(23,2),
\\
(37,1),(37,3),(39,1),(168,2),(168,3),(260,1)
\end{array}
\right\}.
\end{equation}

By the \textsc{Magma} function \texttt{HadamardDatabase}~\cite{Magma},
one can construct $219$ inequivalent Hadamard matrices of order $36$,
including the ten Hadamard matrices in~\cite{T86}.
We verified that there is no  pair of inequivalent Hadamard matrices among 
the $89$ inequivalent Hadamard matrices  $H_{F,i,j}$
given in Proposition~\ref{prop:2}, 
the $13$ inequivalent Hadamard matrices $H_{F,s,t}^T$ for 
$(s,t) \in \mathcal{S}_{s,t}$ in~\eqref{eq:T}
and the above $219$ inequivalent Hadamard matrices in~\cite{Magma}.

\subsection{Related codes}

For the $89$ inequivalent Hadamard matrices  $H_{F,i,j}$
given in Proposition~\ref{prop:2},
we calculated the dimensions $\dim(C(H_{F,i,j}))$ and the minimum
weights $d(C(H_{F,i,j}))$ of their codes $C(H_{F,i,j})$.  The results are as follows:
\begin{align*}
\dim(C(H_{F,i,j})) &=
\begin{cases}
16 & \text{ if }(i,j)=(11,1),
(27,1),
(37,3),
(38,2),
\\ & \qquad
(39,1),
(49,2),
(56,2),
(64,1),
(168,2),
\\
18 & \text{ otherwise,}
\end{cases}
\\
d(C(H_{F,i,j})) &=
\begin{cases}
12 & \text{ if } (i,j)=(11,1), (38,2), \\
9 & \text{ otherwise.}
\end{cases}
\end{align*}
In addition, we verified  that
\[
d(C(H_{F,s,t}^T)) =
\begin{cases}
12 & \text{ if } (s,t)=(260,1), \\
9 & \text{ otherwise,}
\end{cases}
\]
along with $C(H_{F,260,1}^T) \cong P_{36}$.
Since $|\Aut(H_{F,260,1})|=72$ (see Table~\ref{Tab:Aut}), we have the following:

\begin{prop}
Let $H'$ be the second  Hadamard matrix formed by codewords of weight $36$ in
the Pless symmetry code $P_{36}$,
which is given in~\cite[Theorem~2.4 (i)]{T}.
Then the Hadamard matrix $H_{F,260,1}^T$ is equivalent to  $H'$.
\end{prop}

The above  proposition indicates a relationship between
the two codes $P_{36}$ and $C(H_{F,260,1})$
via the Hadamard matrix $H'$ given in~\cite[Theorem~2.4 (i)]{T}.

We end this note by giving
some observation on the code  $C(H_{F,260,1})$.
Let $W(C(H_{F,260,1}))$ be the set of $144$ codewords of weight $36$ in $C(H_{F,260,1})$.
Since 
\[
\cT(C(H_{F,260,1}))
=(36,36,1,1,1)
\]
(see Table~\ref{Tab:tuple2}), 
it follows that there are $h_i$ and $h'_i \ (i=1,\dots,36)$ such that
\[
W(C(H_{F,260,1}))=
\{h_{i}, 2h_{i}, h'_{i}, 2h'_{i} \mid i=1,2,\ldots,36\}
\]
and the following two $36 \times 36$ matrices
\[
\left(
\begin{array}{ccc}
&\overline{h_{1}}&\\
&\vdots &\\
&\overline{h_{36}}& \\
\end{array}
\right) \text{ and }
\left(
\begin{array}{ccc}
&\overline{h'_{1}}& \\
&\vdots &\\
&\overline{h'_{36}} &\\
\end{array}
\right)
\]
are  Hadamard matrices (see~\eqref{eq:rho} for the notation $\overline{h_i}$).
This means that  the $144$ codewords of weight $36$ in $C(H_{F,260,1})$ correspond to
the rows and their negatives of two Hadamard matrices.
It is worthwhile to determine whether there is a near-extremal self-dual code with
$A_9=816$ satisfying the condition that the $72$ codewords of weight $36$
correspond to 
the rows and their negatives of a Hadamard matrix.

\bigskip
\noindent
\textbf{Acknowledgments.}
This work was supported by JSPS KAKENHI Grant Number 19H01802.


\begin{landscape}

\begin{table}[th]
\caption{Four-negacirculant codes $F_{36,i}$ $(i=1,2,\ldots,132)$}
\label{Tab:4nega1}
\centering
\medskip
{\footnotesize
\begin{tabular}{c|c|c||c|c|c||c|c|c||c|c|c}
\noalign{\hrule height1pt}
$i$ &  $r_A$ &$r_B$ &$i$ &  $r_A$ &$r_B$ 
& $i$ &  $r_A$ &$r_B$ & $i$ &  $r_A$ &$r_B$ \\
\hline
  1& 100211222&012000012& 34& 120121102&111000001& 67& 112112202&201210201&100& 110002201&121022020\\
  2& 100121012&221100001& 35& 101120120&112200010& 68& 122000102&122111221&101& 100220221&020221020\\
  3& 110200102&221100021& 36& 122002022&022200012& 69& 112012122&102200000&102& 101102222&202210221\\
  4& 120201102&111000011& 37& 120012020&111000000& 70& 120220222&001200022&103& 102000222&211022002\\
  5& 121202010&222200100& 38& 110020121&202100021& 71& 112212022&122012010&104& 100210012&111200011\\
  6& 111201011&201200100& 39& 111100210&122200010& 72& 110202212&201211201&105& 100212122&210100020\\
  7& 112010202&012200021& 40& 122002002&122200021& 73& 112121021&200222210&106& 100110210&212220010\\
  8& 102101012&111100002& 41& 121002022&220100120& 74& 121011111&021110110&107& 101210102&212221011\\
  9& 122202021&020200011& 42& 102222212&121100110& 75& 121220221&010011121&108& 100002222&222201002\\
 10& 121111212&122200010& 43& 120122122&212200102& 76& 111122222&121220111&109& 110012102&121221120\\
 11& 102102221&021200010& 44& 121201221&111200110& 77& 111221001&101221210&110& 112010221&011001010\\
 12& 112211120&001200010& 45& 121021012&222100121& 78& 101012100&121112000&111& 102211212&201110210\\
 13& 122012122&222100022& 46& 100122122&101100001& 79& 122120022&012022121&112& 110102222&211011120\\
 14& 121022021&021100001& 47& 121011120&110200020& 80& 110001101&011222020&113& 102021202&201002102\\
 15& 110010221&021200021& 48& 111020210&122100010& 81& 120110222&221110110&114& 112101122&220212020\\
 16& 101001011&121200011& 49& 112202020&121100010& 82& 120021122&120220121&115& 122012221&120021101\\
 17& 121012112&122200120& 50& 121020201&212200001& 83& 112210212&122102100&116& 111120212&001222201\\
 18& 121122211&222100100& 51& 102201002&221200011& 84& 122120211&101122100&117& 100001211&110120210\\
 19& 110211012&212200111& 52& 122012221&121200110& 85& 101001221&012110020&118& 110122002&200120210\\
 20& 122102212&121100101& 53& 120212112&112200101& 86& 111212210&112020110&119& 122212122&010202202\\
 21& 122222022&122100101& 54& 122020202&220200110& 87& 102210120&201200210&120& 112100001&111012200\\
 22& 120220112&121100121& 55& 101120122&101200100& 88& 122002201&100101120&121& 102211111&200010100\\
 23& 100211101&212200001& 56& 122020211&121200111& 89& 121122212&021202001&122& 122201011&101200010\\
 24& 120021202&212200020& 57& 121010122&111200121& 90& 101112222&101101202&123& 121111211&221121110\\
 25& 102201101&112100020& 58& 121020212&221100112& 91& 122100221&100222000&124& 122021021&002002210\\
 26& 122001012&212200001& 59& 122011022&121100112& 92& 110101202&210022100&125& 122222002&121110210\\
 27& 110222022&102200020& 60& 110121021&122100111& 93& 120220211&200010101&126& 100110202&021021120\\
 28& 121210101&101200010& 61& 102121202&001200011& 94& 112111122&000212021&127& 110221220&020020210\\
 29& 102001112&110200012& 62& 100102220&110000010& 95& 110102001&012221001&128& 122011121&202121100\\
 30& 111022122&212100022& 63& 100202222&221222120& 96& 102020022&112021001&129& 111001001&021222100\\
 31& 120120021&222100020& 64& 121201202&012121210& 97& 122222002&100010101&130& 100201202&120020111\\
 32& 112220201&202100010& 65& 110112122&010120221& 98& 120101210&110001110&131& 102110221&020002201\\
 33& 120210012&222100001& 66& 111201022&100122211& 99& 121122010&121112100&132& 102010022&120010211\\
\noalign{\hrule height1pt}
\end{tabular}
}
\end{table}

\begin{table}[th]
\caption{Four-negacirculant codes $F_{36,i}$ $(i=133,134,\ldots,260)$}
\label{Tab:4nega2}
\centering
\medskip
{\footnotesize
\begin{tabular}{c|c|c||c|c|c||c|c|c||c|c|c}
\noalign{\hrule height1pt}
$i$ &  $r_A$ &$r_B$ &$i$ &  $r_A$ &$r_B$ 
& $i$ &  $r_A$ &$r_B$ & $i$ &  $r_A$ &$r_B$ \\
\hline
133& 111022202&211011021&165& 112210121&002000101&197& 112122011&212200011&229& 102011021&211202111\\
134& 102111101&101000101&166& 100221122&111002121&198& 102010111&221200010&230& 100000011&120110200\\
135& 120022101&201022200&167& 120121011&221212020&199& 111000202&101011110&231& 121212012&001112120\\
136& 121120221&001002010&168& 100002222&020121102&200& 120012221&021111220&232& 101212002&122210211\\
137& 112010112&111000020&169& 101011212&022120111&201& 112011201&122022210&233& 100212202&011110020\\
138& 121121222&002022110&170& 100020221&221200011&202& 102011222&102220122&234& 100211021&112102211\\
139& 121020101&011202010&171& 120010222&211202221&203& 100000202&022120100&235& 112022101&211111001\\
140& 101002122&222201000&172& 102121100&001200000&204& 111202101&221021101&236& 100200021&002022200\\
141& 121220021&221102110&173& 122110202&212222001&205& 112111022&010010200&237& 110012021&011200210\\
142& 122020021&221122120&174& 121122121&002211100&206& 110022112&200020110&238& 100122002&112200021\\
143& 122020112&202001200&175& 120122212&102101110&207& 100210111&111011000&239& 120001212&022220001\\
144& 100100202&101012000&176& 120200202&010200010&208& 100110012&220022202&240& 102111011&012112120\\
145& 101222001&011202010&177& 102112100&221110000&209& 112110002&002012210&241& 101022220&210212000\\
146& 100000212&001200120&178& 112110010&102110100&210& 120110022&210200120&242& 112001122&012110111\\
147& 122022012&121110120&179& 101210201&022112000&211& 100201111&011011200&243& 101211112&001221120\\
148& 110100222&111220000&180& 101011112&112202202&212& 112120021&010020110&244& 102102002&112200120\\
149& 120221122&101120210&181& 122212022&002121201&213& 101100022&020200200&245& 121202121&010111201\\
150& 112002101&220102200&182& 112010012&012221111&214& 101100111&221100010&246& 102111211&002200010\\
151& 121210002&220021020&183& 120212002&112011121&215& 102001021&121200210&247& 102112212&201200121\\
152& 122222222&122102212&184& 101010001&021012000&216& 112201212&102210220&248& 100221202&100200121\\
153& 121021011&220122220&185& 102220202&002011021&217& 120020221&010111020&249& 121122222&111201111\\
154& 100001012&002200210&186& 120211010&120111000&218& 100200222&002212201&250& 122221122&111100020\\
155& 121112011&022000001&187& 110202212&012121220&219& 110210021&200202011&251& 102021221&012112101\\
156& 100212102&200022120&188& 111021022&122002000&220& 122110102&002002101&252& 120211202&021200010\\
157& 112201022&200001201&189& 120011011&121000120&221& 111012120&212002000&253& 111100010&021100000\\
158& 111100012&022012010&190& 110220021&022200110&222& 100210202&101002111&254& 112200020&111000000\\
159& 102102020&010011000&191& 122100002&112201200&223& 120200122&021210100&255& 122200010&222000000\\
160& 112121022&001110111&192& 112210012&100120002&224& 120121120&102221110&256& 101010120&110000010\\
161& 100202122&222011121&193& 111201102&110200001&225& 100021112&101000211&257& 101112122&212221000\\
162& 110221112&200211021&194& 122001120&010100000&226& 110020002&000101011&258& 120012001&001212110\\
163& 101001202&120111001&195& 101211222&120220201&227& 110110221&111100221&259& 112210201&220200020\\
164& 102210011&002001000&196& 102112001&121200001&228& 102120212&212121200&260& 112101021&200000000\\
\noalign{\hrule height1pt}
\end{tabular}
}
\end{table}

\begin{table}[p]
\caption{Tuples $\cT(F_{36,i})$}
\label{Tab:tuple}
\centering
\medskip
{\footnotesize
\begin{tabular}{c|c||c|c||c|c}
\noalign{\hrule height1pt}
\multicolumn{1}{c|}{$i$}  & \multicolumn{1}{c||}{$\cT(F_{36,i})$}  &
\multicolumn{1}{c|}{$i$}  & \multicolumn{1}{c||}{$\cT(F_{36,i})$}  &
\multicolumn{1}{c|}{$i$}  & \multicolumn{1}{c}{$\cT(F_{36,i})$}  \\
\hline

$1,6$  &$(144,264,0,0,0)$ &$35$   &$(120,216,0,1,1)$ &$80,81,87,98,99$&$(264,36,0,0,0)$\\
$2,5$  &$(408,0,0,0,0)$   &$37$   &$(84,252,0,19,3)$ &$83$&$(216,84,1,0,1)$\\
$3$    &$(372,36,0,0,0)$  &$38$   &$(72,264,1,1,2)$  &$88,94$&$(48,252,0,0,0)$\\
$4$    &$(264,144,1,1,2)$ &$39$   &$(288,48,3,0,1)$  &$89$&$(36,264,0,0,0)$\\
$7$    &$(252,144,3,0,1)$ &$40,44$&$(108,228,0,0,0)$ &$91$&$(156,144,0,0,0)$\\
$8$    &$(288,108,12,0,2)$&$41$   &$(144,192,0,0,0)$ &$95,101$&$(84,216,0,0,0)$\\
$9$    &$(360,12,1,0,1)$  &$42$   &$(180,156,1,0,1)$ &$103$&$(84,216,0,1,1)$\\
$10$   &$(144,228,0,0,0)$ &$43$   &$(264,72,0,0,0)$  &$105,109,113$&$(108,180,0,0,0)$\\
$11$   &$(288,84,1,0,1)$  &$45$   &$(72,264,0,0,0)$  &$106,107,116$&$(216,72,0,0,0)$\\
$12$   &$(264,108,4,0,2)$ &$46$   &$(144,180,0,0,0)$ &$108$&$(216,72,2,0,2)$\\
$13$   &$(84,288,0,1,1)$  &$47$&$(252,72,0,0,0)$&$110$&$(216,72,1,0,1)$\\
$14$   &$(300,72,0,0,0)$  &$48,52$&$(288,36,0,0,0)$&$111$&$(180,108,1,0,1)$\\
$15$   &$(336,36,0,0,0)$  &$49$&$(144,180,0,2,2)$&$112,117,120$&$(36,252,0,1,1)$\\
$16$   &$(72,300,0,0,0)$  &$50,57$&$(180,144,0,0,0)$&$114$&$(180,108,0,0,0)$\\
$17$   &$(216,156,0,0,0)$ &$51,55,58,59$&$(216,108,0,0,0)$&$115$&$(108,180,0,1,1)$\\
$18$   &$(84,288,0,14,3)$ &$53$&$(36,288,0,1,1)$&$118,119$&$(144,144,0,0,0)$\\
$19,21$&$(36,336,0,0,0)$  &$54$&$(36,288,0,2,2)$&$121$&$(72,216,0,0,0)$\\
$20$   &$(228,144,0,0,0)$ &$56$&$(252,72,2,0,2)$&$122$&$(192,72,0,0,0)$\\
$22$   &$(192,180,0,0,0)$ &$60$&$(288,36,1,0,1)$&$123,129,132,137,141,154$&$(144,120,0,0,0)$\\
$23$   &$(252,108,19,0,3)$&$61,66$&$(264,36,1,0,1)$&$124,138$&$(36,228,0,0,0)$\\
$24$   &$(216,144,0,0,0)$ &$62,69,71,75,77$&$(144,156,0,0,0)$&$125,128,144,146$&$(84,180,0,0,0)$\\
$25$   &$(72,288,0,0,0)$  &$63,65,72,84,85,97,100$&$(192,108,0,0,0)$&$126,140,152$&$(72,192,0,0,0)$\\
$26$   &$(36,324,0,4,2)$  &$64,90$&$(144,156,0,1,1)$&$127,134,150$&$(108,156,0,0,0)$\\
$27$   &$(156,180,0,1,1)$ &$67,96$&$(156,144,0,1,1)$&$130,139$&$(84,180,0,1,1)$       \\
$28$   &$(300,36,2,0,1)$  &$68$&$(72,228,0,1,1)$&$131,135$&$(120,144,0,0,0)$            \\
$29$   &$(120,216,0,0,0)$ &$70,74,92,93,102$&$(108,192,0,0,0)$&$133$&$(120,144,0,1,1)$    \\
$30$   &$(288,48,1,0,1)$  &$73$&$(48,252,0,1,1)$&$136,147,153$&$(180,84,0,0,0)$             \\
$31,32$&$(156,180,0,0,0)$ &$76,82,86$&$(180,120,0,0,0)$&$142,151$&$(228,36,0,0,0)$            \\
$33,36$&$(216,120,0,0,0)$ &$78$&$(216,84,0,0,0)$&$143$&$(72,192,0,2,2)$                     \\
$34$   &$(48,288,0,0,0)$  &$79,104$&$(192,108,1,0,1)$&$145$&$(156,108,0,0,0)$                 \\
\noalign{\hrule height1pt}
\end{tabular}
}
\end{table}

\begin{table}[htbp]
\caption{Tuples $\cT(F_{36,i})$}
\label{Tab:tuple2}
\centering
\medskip
{\footnotesize
\begin{tabular}{c|c||c|c||c|c}
\noalign{\hrule height1pt}
\multicolumn{1}{c|}{$i$}  & \multicolumn{1}{c||}{$\cT(F_{36,i})$}  &
\multicolumn{1}{c|}{$i$}  & \multicolumn{1}{c||}{$\cT(F_{36,i})$}  &
\multicolumn{1}{c|}{$i$}  & \multicolumn{1}{c}{$\cT(F_{36,i})$}  \\
\hline
$148$&$(228,36,1,0,1)$                     &$191$&$(48,180,0,0,0)$&$230$&$(144,48,1,0,1)$\\
$149$&$(192,72,1,0,1)$                     &$193$&$(192,36,1,0,1)$       &$231$&$(0,192,0,0,0)$\\
$155,157,167$&$(144,108,0,0,0)$            &$194$&$(144,84,1,0,1)$       &$234$&$(84,108,0,1,1)$\\
$156,158,160,161,165,171$&$(108,144,0,0,0)$&$197$&$(36,180,0,1,1)$       &$235,237$&$(108,84,1,0,1)$\\
$159,163$&$(216,36,0,0,0)$                 &$198,199,204,212$&$(108,108,0,0,0)$  &$238,239,241$&$(144,36,0,0,0)$\\
$162,170$&$(72,180,0,0,0)$                 &$200$&$(36,180,0,0,0)$&$240$&$(36,144,0,0,0)$\\
$164$&$(72,180,0,1,1)$                     &$201$&$(108,108,0,1,1)$&$242,244$&$(108,72,0,0,0)$\\
$166$&$(216,36,1,0,1)$                     &$202,205,209$&$(180,36,0,0,0)$&$243$&$(72,108,0,0,0)$\\
$168$&$(252,0,20,0,4)$                     &$203,206,211$&$(144,72,0,0,0)$&$245$&$(144,36,2,0,1)$\\
$169$&$(252,0,1,0,1)$                      &$207$&$(72,144,0,0,0)$&$246,247$&$(120,36,1,0,1)$\\
$172,188$&$(0,228,0,1,1)$                  &$208$&$(108,108,1,0,1)$&$248$&$(108,48,0,0,0)$\\
$173$&$(72,156,0,1,1)$                     &$210$&$(180,36,1,0,1)$&$249$&$(72,84,0,0,0)$\\
$174$&$(228,0,1,0,1)$                      &$213$&$(144,72,1,1,2)$&$250,251$&$(72,72,0,0,0)$\\
$175,180,187,189$&$(156,72,0,0,0)$         &$214,218,219,225,236$&$(108,84,0,0,0)$&$252$&$(72,48,0,0,0)$\\
$176$&$(36,192,0,2,2)$                     &$215,224,226,227$&$(36,156,0,1,1)$&$253$&$(48,72,0,0,0)$\\
$177$&$(156,72,1,0,1)$                     &$216$&$(120,72,0,0,0)$&$254$&$(12,108,0,2,1)$\\
$178$&$(48,180,0,1,1)$                     &$217,233$&$(84,108,0,0,0)$&$255$&$(120,0,0,0,0)$\\
$179,182,185,192,195$&$(108,120,0,0,0)$    &$220,222,228,232$&$(36,156,0,0,0)$&$256$&$(72,12,0,0,0)$\\
$181,183,196$&$(120,108,0,0,0)$            &$221$&$(144,48,0,0,0)$&$257,258$&$(36,48,0,0,0)$\\
$184$&$(108,120,0,1,1)$                    &$223$&$(48,144,0,0,0)$&$259$&$(72,0,0,0,0)$\\
$186,190$&$(72,156,0,0,0)$                 &$229$&$(72,120,0,0,0)$&$260$&$(36,36,1,1,1)$\\
\noalign{\hrule height1pt}
\end{tabular}
}
\end{table}

\end{landscape}

\end{document}